\documentclass[12pt]{article}

\usepackage{graphicx}
\usepackage{epsfig}
\usepackage{epstopdf}

\usepackage{amssymb, amsmath, amsthm}

\newtheorem{dfn}{Definition}
\newtheorem{defn}[dfn]{Definition}

\newtheorem{rem}[dfn]{Remark}

\newtheorem{thm}{Theorem}
\newtheorem{lem}{Lemma}

\newtheorem{prop}[thm]{Proposition}

\newtheorem{ex}{Example}

\def\proof{\par\medskip\noindent{\it Proof: }}

\def\>{\rangle}
\def\<{\langle}

\def\3{\ss}
\def\8{\infty}

\begin{document}


\title{Kakimizu complexes of Surfaces and 3-Manifolds}
\author{Jennifer Schultens}

\maketitle

\begin{abstract}
The Kakimizu complex is usually defined in the context of
knots, where it is known to be quasi-Euclidean.   We here generalize the definition of the Kakimizu complex to surfaces and 3-manifolds 
(with or without boundary).  Interestingly, in the setting of surfaces, the complexes and the techniques turn out to replicate those used to study 
the Torelli group, {\it i.e.,} the ``nonlinear" subgroup of the mapping class group.  Our main results are that the Kakimizu complexes of a surface are contractible and that they need not be quasi-Euclidean.
It follows that there exist (product) $3$-manifolds whose Kakimizu complexes are not quasi-Euclidean.  


\end{abstract}

\vspace{2 mm}

The existence of Seifert's algorithm, discovered by Herbert Seifert, proves, among other things, that every 
knot admits a Seifert surface.  {\it I.e.,} for every knot $K$, there is a compact orientable surface whose boundary
is $K$.  It is worth noting that the existence of a Seifert surface for a knot $K$ also follows 
from the existence of submanifolds representing homology classes of manifolds or pairs of submanifolds, in this 
case the pair $(K, {\mathbb S}^3)$.   This point of view proves useful in generalizing our understanding of
Seifert surfaces to other classes of surfaces in $3$-manifolds.  

Adding a trivial handle to a Seifert surface produces an isotopically distinct 
surface.  Adding additional handles produces infinitely many isotopically distinct surfaces.
These are not the multitudes of surfaces of primary interest here.
The multitudes of surfaces of primary interest here are, for example, the infinite collection of Seifert surfaces produced by
Eisner, see \cite{E}.  Eisner realized that ``spinning" a Seifert surface around the decomposing annulus of a
connected sum of two non fibered knots produces homeomorphic but non isotopic Seifert surfaces. 
This abundance of Seifert surfaces led Kakimizu to define a complex, now named after him, whose vertices
are isotopy classes of Seifert surfaces of a given knot and whose $n$-simplices are $(n+1)$-tuples of vertices that admit 
pairwise disjoint representatives.  


Our understanding of the topology and geometry of the Kakimizu complex continues to evolve.  Both
Kakimizu's work and, independently, a result of Scharlemann and Thompson, imply that the Kakimizu complex is
connected.  (See \cite{K} and \cite{ST}.)  Sakuma and Shackleton exhibit diameter bounds in terms of the genus of a knot.  
(See \cite{SS}.)  P. Przytycki and the author establish that the Kakimizu complex is contractible.  (See \cite{PS}.)
Finally, Johnson, Pelayo and Wilson prove that the Kakimizu complex of a knot is quasi-Euclidean.  (See \cite{JPW}.)

This paper grew out of a desire to study concrete examples of Kakimizu complexes of $3$-manifolds other than 
knot complements.  A natural case to consider is product manifolds, where relevant information is captured
by the surface factor.  The challenge lies in adapting the idea of the Kakimizu complex to a more general setting: 
codimension $1$ submanifolds of $n$-manifolds.   

As it turns out, in the case of $1$-dimensional submanifolds of a surface, the Kakimizu complexes are related to 
the homology curve complexes investigated by Hatcher (see \cite{Hatcher}), Irmer (see \cite{II}), 
Bestvina-Bux-Margalit (see \cite{BBM}) and Hatcher-Margalit (see \cite{HM}) discussed in Section 2.  
These complexes are of interest in the study of the Torelli group, which is the kernel of the action of the mapping class group 
of a manifold on the homology of the manifold.    The Torelli group of a surface, in turn, acts on the homology curve complexes. 
This group action has been used to study the topology of the Torelli group of a surface, 
for instance by Bestvina-Bux-Margalit in their investigation of the dimension of the 
Torelli group (see \cite{BBM}), by Irmer in ``The Chillingworth class is a signed stable length'' (see \cite{II-}), 
by Hatcher-Margalit in ``Generating the Torelli group'' (see \cite{HM}), and by Putman in 
``Small generating sets for the Torelli group'' (see \cite{P}).  


Hatcher proved that the homology curve complex is contractible and computed its dimension.  Irmer studied geodesics
of the homology curve complex and exhibited quasi-flats.  These insights guide our investigation of the 
Kakimizu complex of a surface.   Specifically, we prove similar, and in some cases analogous, results in the setting
of the Kakimizu complex of a surface.   Our main results are that the Kakimizu complexes of a surface are contractible and 
that they need not be quasi-Euclidean.   

One example stands out:  The Kakimizu complex of a genus $2$ surface. 
In \cite{BBM}, Bestvina-Bux-Margalit reprove a theorem of Mess, that the Torelli group of a genus $2$ surface is an infinitely
generated free group.  They do so by showing that it acts on a tree with infinitely many edges emerging from each vertex.
As it turns out, the Kakimizu complex of the genus $2$ surface is also a tree with infinitely many edges emerging from
each vertex.  In particular, the Kakimizu complex of the genus $2$ surface is Gromov hyperbolic.   A product manifold
with the genus $2$ surface as a factor will thus also have some Gromov hyperbolic Kakimizu complexes.  This is interesting as it
shows that in addition to examples of $3$-manifolds with quasi-Euclidean Kakimizu complexes, as proved by Johnson-Pelayo-Wilson,
there are  $3$-manifolds with Gromov hyperbolic Kakimizu complexes.  Kakimizu complexes exhibit more than one geometry!


I wish to thank Misha Kapovich, Allen Hatcher, Dan Margalit and Piotr Przytycki for helpful conversations and the referee for helpful comments.  
This research was supported by a grant from the NSF.  Much of the work was carried out at the 
Max Planck Institute for Mathematics in Bonn.  I thank the Max Planck Institute for Mathematics in Bonn for its hospitality
and the Hausdorff Institute for Mathematics for logistical support.   



\section{The Kakimizu complex of a surface}

The work here follows in the footsteps of \cite{PS}.  
Whereas the setting for \cite{PS} is surfaces in $3$-manifolds, the setting here is $1$-manifolds in $2$-manifolds.
It is worth pointing out that although we discuss only $1$-manifolds in $2$-manifolds and $2$-manifolds in $3$-manifolds, 
the definitions and arguments carry over verbatim to the setting of codimension $1$ submanifolds in manifolds of any dimension.

Recall that an element of a finitely generated free abelian group $G$ is
{\em primitive} if it is an element of a basis for $G$.  
In the following we will always assume: 
1) $S$ is a compact (possibly closed) connected oriented 2-manifold; 
2) $\alpha$ is a primitive element of $H_1(S, \partial S, {\bf Z})$.  


\begin{defn} \label{Seifert-curve}
A {\em Seifert curve} for $(S, \alpha)$ is a pair $(w, c)$, where $c$ is a union, $c_1 \sqcup \cdots \sqcup c_n$, of
pairwise disjoint oriented 
simple closed curves and arcs in $S$ and $w$ is an $n$-tuple of natural numbers $(w^1, \dots, w^n)$ such
that the homology class $w^1[[c_1]] + \cdots + w^n[[c_n]]$ equals $\alpha$.   
Moreover, we require that $S \backslash c$ is connected.  We call $c$ the underlying curve of $(w , c)$.  
We will denote $w^1[[c_1]] + \cdots + w^n[[c_n]]$ by $w \circ c$.  
\end{defn}

\begin{figure}[ht]
\vspace{2 mm}
\centerline{\epsfxsize=4in \epsfbox{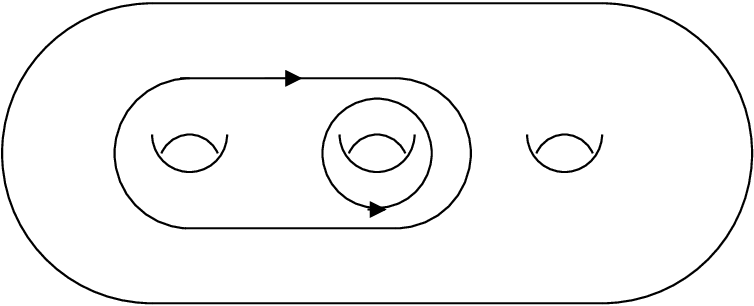}}
\caption{\sl A Seifert curve (weights are $1$)}
\label{diskarg}
\vspace{2 mm}
\end{figure}

Our definition of Seifert curve disallows null homologous subsets. Indeed, a null homologous subset would bound a component of 
$S \backslash c$ and would hence be separating.  In fact, $c$ contains no bounding subsets.
Conversely, if $w \circ d = \alpha$ and $d$ contains no bounding subsets, then $S \backslash d$ is connected.  

\begin{lem} \label{uniquem}
If $(w, c)$ represents $\alpha$, then $w$ is determined by the underlying curve $c$. 
\end{lem}

\proof
Suppose that $(w, c)$ and $(w', c)$ represent $\alpha$, where $w = (w^1, \cdots, w^n)$ and
$w' = ((w')^1, \cdots, (w')^n)$.  Then \[w^1[[c_1]] + \cdots + w^n [[c_n]] = \alpha = (w')^1[[c_1]] + \cdots + (w')^n [[c_n]],\]
hence \[(w^1-(w')^1)[[c_1]] + \cdots (w^n - (w')^n)[[c_n]] = 0.\]  Since $c$ has no null homologous subsets,
this ensures that \[w^1 - (w')^1 =  0, \dots, w^n-(w')^n = 0.\]
Thus \[w^1 =  (w')^1, \cdots, w^n = (w')^n.\]
\qed

Since the 
underlying curve $c$ of a Seifert curve $(w, c)$ determines $w$, we will occasionally speak of a Seifert curve
$c$, when $w$ does not feature in our discussion. 

\begin{defn}
Given a Seifert curve $(w, c)$ we denote the curve obtained by replacing, for all $i$, the curve $c_i$ with 
$w_i$ parallel components of $c_i$, by $h(w, c)$.  This defines a function from Seifert curves to unweighted curves.

Conversely, let $d = d_1 \sqcup \dots \sqcup d_m$ be a disjoint union of (unweighted) pairwise disjoint simple
closed curves and arcs such that parallel components are oriented to be parallel oriented curves and arcs.  
We denote the weighted curve obtained by replacing parallel components with one weighted
component whose weight is equivalent to the number of these parallel components by $h^{-1}(d)$.  
\end{defn}

\begin{defn}
For each pair $(S, \alpha)$, the isomorphism between $H_1(S, \partial S)$ and $H^1(S)$ identifies
an element $a^*$ of $H^1(S)$ corresponding to $\alpha$ that lifts to a homomorphism $h_a: \pi_1(S) \rightarrow {\bf Z}$.
We denote the covering space corresponding to $N_{\alpha} = kernel(h_a)$ by $(p_{\alpha}, \hat S_{\alpha}, S)$, or simply
$(p, \hat S, S)$, and call it the {\em infinite cyclic covering space associated with $\alpha$}.
\end{defn}

\begin{figure}[h]
\vspace{2 mm}
\centerline{\epsfxsize=4in \epsfbox{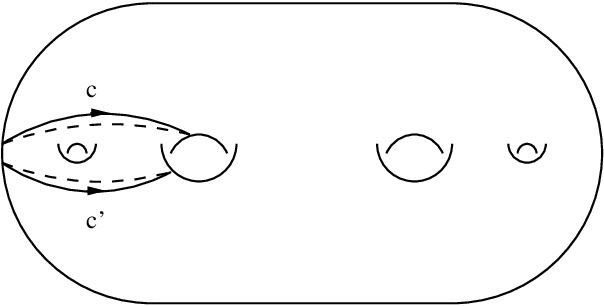}}
\caption{\sl Two Seifert curves corresponding to vertices of distance 1 (weights are 1)}
\label{cc}
\vspace{2 mm}
\end{figure}

We now describe the Kakimizu complex of $(S, \alpha)$.  As {\em vertices}
we take Seifert curves $(w, c)$ of $(S, \alpha)$, considered up to isotopy of underlying curves.  
We write $[(w,c)]$.  
Consider a pair of vertices $v, v'$ and representatives $(w, c), (w', c')$.
Here $S \backslash c$ and $S \backslash c'$ are connected, hence path-connected.
It follows that lifts of $S \backslash c$ and $S \backslash c'$ to the covering space associated with $\alpha$ are simply path components
of $p^{-1}(S \backslash c)$ and $p^{-1}(S \backslash c')$.  
We obtain a graph $\Gamma(S, \alpha)$ by spanning an {\em edge} $e=(v, v')$ on the vertices $v, v'$ if and only if 
the representatives $(w, c), (w',c')$ of $v, v'$ can be chosen 
so that a lift of $S \backslash c$ to the covering space associated with $\alpha$
intersects exactly two lifts of $S \backslash c'$.  (Note that in this case
$c$ and $c'$ are necessarily disjoint.)   See Figure \ref{cc}.  

\begin{defn} 
Let $X$ be a simplicial complex.  If, whenever the $1$-skeleton of a simplex $\sigma$ is in $X$, the simplex $\sigma$ is also 
in $X$, then $X$ is said to be {\em flag}.  
\end{defn}

\begin{defn} \label{Kakimizu-complex}
The {\em Kakimizu complex of $(S, \alpha)$}, denoted by 
$Kak(S, \alpha)$ is the flag complex with $\Gamma(S, \alpha)$ as its $1$-skeleton. 
\end{defn}

\begin{rem}
The Kakimizu complex is defined for a pair $(S, \alpha)$.  For simplicity we use the expression ``the Kakimizu
complex of a surface" in general discussions, rather than the more cumbersome ``the Kakimizu complex of a pair $(S, \alpha)$,
where $S$ is a surface and $\alpha$ is a primitive element of $H_1(S, \partial S, {\bf Z})$".   Note that the Kakimizu complex of a surface 
is thus unique only in conjunction with a specified $\alpha$.  
\end{rem}

\begin{figure}[h]
\vspace{2 mm}
\centerline{\epsfxsize=4in \epsfbox{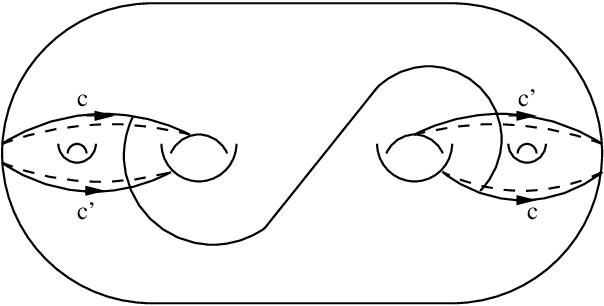}}
\caption{\sl Two Seifert curves corresponding to vertices of distance strictly greater than 1 (weights are 1)}
\label{ccstar}
\vspace{2 mm}
\end{figure}

Figure \ref{ccstar} provides an example of a pair $(w, c), (w', c')$ of disjoint (disconnected) Seifert curves that do not span an edge.  
The arc from one side of $c$ to the other side of $c$ intersects $c'$ twice with the same orientation and a
lift of $S \backslash c$ will 
hence meet at least three distinct lifts of $S \backslash c'$.  For a 3-dimensional analogue of Figure \ref{ccstar}, 
see \cite{Banks2011}.

\begin{ex}
The Kakimizu complexes of the disk and sphere are empty.  The annulus has a non-empty but
trivial Kakimizu complex $Kak(A, \alpha)$ consisting of a single vertex.  Specifically, 
let $A = annulus,$ and $\alpha$ a generator of 
$H_1(A, \partial A, {\mathbb Z}) = {\mathbb Z}$.  Then $\alpha$ is represented by a spanning arc with weight $1$.  
The spanning arc is, up to isotopy, the only possible underlying curve for a representative of $\alpha$.
Thus $Kak(A, \alpha)$ consists of a single vertex.  
 
Similarly, the torus has non-empty but
trivial Kakimizu complexes, each consisting of a single vertex.  Specifically, let $T = torus,$
and $\beta$ a primitive element of $H_1(T, {\mathbb Z}) = {\mathbb Z} \times {\mathbb Z}$.  
Again, there is, up to isotopy, only one underlying curve for representatives of $\beta$.  There are infinitely many choices for $\beta$, but
in each case, $Kak(T, \beta)$ consists of a single vertex.  
\end{ex}

Having understood the above Examples, we restrict our attention to the case where $S$ is a compact orientable hyperbolic surface with geodesic boundary for the remainder of this paper.  

\begin{defn}
Let $(w, c)$ and $(w', c')$ be Seifert curves.  We say that $(w, c)$ and $(w', c')$ (or simply $c$ and $c'$)
are {\em almost disjoint} if for all $i, j$ the component $c_i$ of $c$ and the component $c_j'$ of $c'$ are either disjoint or coincide.  
\end{defn}

\begin{rem} \label{simplex}
Let $\sigma$ be a simplex in $Kak(S, \alpha)$ of dimension $n$.  Denote the vertices of $\sigma$ by 
$v_0, \dots, v_n$ and let $c_0, \dots, c_n$ be geodesic representatives of the underlying curves of 
Seifert curves for 
$v_0, \dots, v_n$ such that 
arc components of $c_0, \dots, c_n$ are perpendicular to $\partial S$.   It is a well known fact that closed
geodesics that can be isotoped to be disjoint must be disjoint or coincide.  The same is true for the geodesic arcs
considered here, and combinations of closed geodesics and geodesic arcs, because their doubles are 
closed geodesics in the double of $S$.   Hence, for all pairs $i, j$, the component $c_i$ of $c$ and the component
$c_j'$ of $c'$ are either disjoint or coincide.
\end{rem}

\begin{defn}  \label{d_K}
Consider $Kak(S, \alpha)$.  Let $(p, \hat S, S)$ be the infinite cyclic cover of $S$ 
associated with $\alpha$.
Let $\tau$ be the generator of the group of 
covering transformations of $(p, \hat S, S)$ (which is ${\mathbb Z}$) corresponding to $1$.  
Note that $\tau$ is canonical up to sign.  

Let $(w, c), (w', c')$ be Seifert curves in 
$(S, \alpha)$.  
Let $S_0$ denote a lift of $S \backslash c$ to 
$\hat S$, {\it i.e.}, a path component of $p^{-1}(S \backslash c)$.
Set $S_i = \tau^i(S_0)$,
$c_i =  closure(S_i) \cap closure(S_{i+1})$.  
Let $S_0'$ be a lift of 
$S \backslash c'$ to $\hat S$.  Set $d_K(c, c) = 0$ and for $c \neq c'$, set 
{\em $d_K(c, c')$} equal to one less than the number of translates of $S_0$ met by $S_0'$.  
Let $v, v'$ be vertices in $Kak(S, \alpha)$.  Set $d_k(v, v) = 0$ 
and for $v \neq v'$ set {\em $d_K(v, v')$} equal to the minimum of 
$d_K(c, c')$ for $(w, c), (w', c')$ representatives of $v, v'$.  
\end{defn}

\begin{defn} \label{belowabove}
Let $C, D$ be disjoint separating subsets of $\hat S$.
We say that $D$ lies {\em above}
$C$ if $D$ lies in the component of $\hat S \backslash C$ containing $\tau(C)$.  We
say that $D$ lies {\em below} $C$ if $D$ lies in the component of $\hat S \backslash C$ containing $\tau^{-1}(C)$.   
\end{defn}

\begin{rem} \label{d_Kisfinite}
Here $d_K(c, c')$ is finite:  Indeed, $w \circ c = w' \circ c' = \alpha$ and so $[(w, c)],$ 
$[(w', c')]$ are in the kernel, $N_{\alpha}$, of the cohomology class dual to $\alpha$.  
Specifically, the 
cohomology class dual to $\alpha$ is represented
by the weighted intersection pairing with $(w, c)$ and also the weighted 
intersection pairing with $(w', c')$.  Thus, let $c^j$ be a component of $c$, then
the value of the cohomology class dual to $\alpha$ evaluated at $[[c^j]]$ is given by 
the weighted intersection pairing of $(w, c)$ with $c^j$
which is $0$.  Likewise for other components of $c$ and $c'$.  
Thus each component of $c, c'$ lies in the kernel of this homomorphism and hence 
in $N_{\alpha}$.  Thus lifts of $c, c'$ are homeomorphic to $c, c'$, respectively,
in particular, they are compact 1-manifolds.   It follows that $d_K(c, c')$ is finite, whence
for all vertices $v, v'$ of $Kak(S, \alpha)$, $d_K(v, v')$ is also finite.
\end{rem}

It is not hard to verify, but important to note, the following theorem (see \cite[Proposition 1.4]{K}):

\begin{prop}
The function $d_K$ is a metric on the vertex set of $Kak(S, \alpha)$.  
\end{prop}


\section{Relation to homology curve complexes} \label{compare}


In \cite{Hatcher}, Hatcher introduces the {\em cycle complex} of a surface:

``By a {\em cycle} in a closed oriented surface $S$ we mean a nonempty collection of finitely many 
disjoint oriented smooth simple closed curves. A cycle $c$ is {\em reduced} if no subcycle of $c$ is 
the oriented boundary of one of the complementary regions of $c$ in $S$ (using either orientation of the region). 
In particular, a reduced cycle contains no curves that bound disks in $S$, and no pairs of circles that are parallel 
but oppositely oriented.

Define the {\em cycle complex} $C(S)$ to be the simplicial complex having as its vertices the isotopy classes of 
reduced cycles in $S$, where a set of $k + 1$ distinct vertices spans a $k$-simplex if these vertices are represented 
by disjoint cycles $c_0,$ $\dots,$ $c_k$ that cut $S$ into $k + 1$ cobordisms $C_0,$ $\dots,$ $C_k$ 
such that the oriented boundary of $C_i$ is $c_{i+1}-c_i$, subscripts being taken modulo $k + 1$, where the orientation 
of $C_i$ is induced from the given orientation of $S$ and $-c_i$ denotes $c_i$ with the opposite orientation. 
The cobordisms $C_i$ need not be connected. The faces of a $k$-simplex are obtained by deleting a cycle 
and combining the two adjacent cobordisms into a single cobordism. One can think of a $k$-simplex of $C(S)$ 
as a cycle of cycles. The ordering of the cycles $c_0,$ $\dots,$ $c_k$ in a $k$-simplex is determined 
up to cyclic permutation.  Cycles that span a simplex represent the same element of $H_1(S)$ since they are 
cobordant. Thus we have a well-defined map $\pi_0: C(S) \rightarrow H_1(S)$. 
This has image the nonzero elements of $H_1(S)$ since on the one hand, every cycle representing 
a nonzero homology class contains a reduced subcycle representing the same class 
(subcycles of the type excluded by the definition of reduced can be discarded one by one until 
a reduced subcycle remains), and on the other hand, it is an elementary fact, left as an exercise, 
that a cycle that represents zero in $H_1(S)$ is not reduced.
For a nonzero class $x \in H_1(S)$ let $C_x(S)$ be the subcomplex of $C(S)$ spanned by
vertices representing $x$, so $C_x(S)$ is a union of components of $C(S)$."  See \cite[Page 1]{Hatcher}.  

\begin{lem} \label{subcomplex}
When both are defined, {\it i.e.,} when $S$ is closed, connected, of genus at least $2$, and $\alpha$ is primitive, 
$Vert(Kak(S, \alpha))$ is isomorphic to a proper subset of $Vert(C_{\alpha}(S))$.  
\end{lem}

\proof
Let $v$ be a vertex of $Kak(S, \alpha)$.  If we choose a representative $(w, c)$, then $h(w, c)$ is
a disjoint collection of (unweighted) curves and arcs.  The requirement on the Seifert curve $(w, c)$, 
that $S \backslash c$ be connected implies that the multi-curve $h(w, c)$ is reduced and thus
represents a vertex of $C_{\alpha}(S)$.  Abusing notation slightly, we denote the map
from $Vert(Kak(S, \alpha))$ to $Vert(C_{\alpha})$ thus obtained by $h$.  There is an inverse,
$h^{-1}$, defined on the image of $h$, hence $h$ is injective. 

It is not hard to identify reduced multi-curves that contain bounding subsets that are not the oriented 
boundary of a subsurface.  Hence $Vert(Kak(S, \alpha))$ is a proper subset of $Vert(C_{\alpha}(S))$.  
\qed

\begin{lem} \label{cycleofcycles}
Suppose that $S$ is hyperbolic and let $\sigma$ be an $n$-simplex in $Kak(S, \alpha)$.  Denote the vertices of $\sigma$ by $v_0, \dots, v_n$.  Then there are representatives of $v_0, \dots, v_n$ with underlying curves
$c_0, \dots, c_n$ such that the following hold:

\begin{enumerate}
\item $c_i \cap c_j = \emptyset$ $\forall i \neq j$; 
\item $S \backslash (c_0 \cup \dots \cup c_n)$ is partitioned into subsurfaces 
$P_0, \dots, P_n$ such that
$\partial P_i = c_i  -  c_{i-1}$.
\end{enumerate}
\end{lem}

\proof
Let $(p, \hat S, S)$ be the covering space associated with $\alpha$ and let $\sigma$ 
be a simplex in $Kak(S, \alpha)$.   Let $c_0, \dots, c_n$ be geodesic representatives of 
the underlying curves of $v_0, \dots, v_n$ such that 
arc components of $c_0, \dots, c_n$ are perpendicular to $\partial S$.  By Remark \ref{simplex},
$c_i$ and $c_j$ are almost disjoint $\forall i \neq j$.  Consider a lift $S_0$ of $S \backslash c_0$ 
to $\hat S$.  For each $j \neq 0$, $c_j$ lifts to a separating collection $\hat c_j$ of simple closed 
curves and simple arcs.  Moreover, since $S_0$ is homeomorphic to $S \backslash c_0$, the lifts 
$\hat c_i, \hat c_j$ are
almost disjoint as long as $i \neq j$.  
By reindexing $c_0, \dots, c_n$ if necessary and performing small isotopies that pull apart equal components, 
we can thus ensure that
$\hat c_i$ lies above $\hat c_j$ for $i > j$.  

Note that the lift of $S \backslash c_0$ is homeomorphic to $S \backslash c_0$.
In particular, $c_i \cap c_j = \emptyset$ $\forall i \neq j$.
Moreover, the surface with interior below $\hat c_i$ and above $\hat c_{i-1}$ projects to a subsurface $P_i$ of $S$ for $i = 1, \dots, n$.
The subsurfaces $P_1, \dots, P_n$ 
exhibit the required properties.
\qed

\begin{figure}[h]
\vspace{2 mm}
\centerline{\epsfxsize=3in \epsfbox{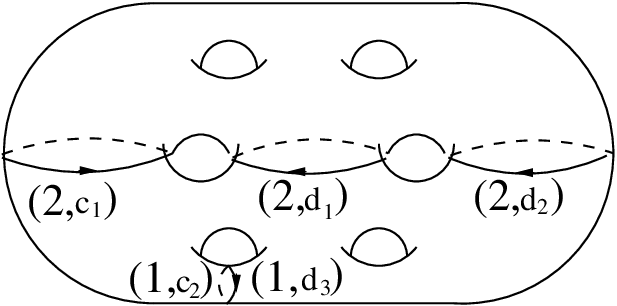}}
\caption{\sl An edge of $Kak(S, \alpha)$ that does not map into $C_{\alpha}$}
\label{nothatcher}
\vspace{2 mm}
\end{figure}

\begin{rem}
When $w_0, \dots, w_n = 1$, Lemma \ref{cycleofcycles} ensures that $c_0$ $=$ $h(w_0, c_0)$ $=$ $h(1, c_0),$ $\dots,$ $c_n$ $=$ $h(w_n, c_n)$ $=$ $h(1, c_n)$ 
form a cycle of cycles.  In this case $h$ extends over the simplex $\sigma$ to produce a simplex 
$h(\sigma)$ in $C_{\alpha}$.  
However, $h$ does not extend over simplices in which weights are not all $1$.  See Figure
\ref{nothatcher}.  
\end{rem}

Hatcher proves that for each $x \in H_1(S)$, $C_x(S)$ is contractible.  (In particular, it is therefore
connected and hence constitutes just one component of $C(S)$.) In Section \ref{contractibility} we
prove an analogous result for $Kak(S, \alpha)$, using a technique from the study of the Kakimizu complex of 3-manifolds.

The cyclic cycle complex and the Kakimizu complex are simplicial complexes.  The complex defined by
Bestvina-Bux-Margalit (see \cite{BBM}) is not simplicial, but can be subdivided to obtain a simplicial complex.  
See the final comments in \cite[Section 2]{HM}.  There is a subcomplex of the cyclic cycle complex that equals
this subdivision of the complex defined by Bestvina-Bux-Margalit.   This is the complex of interest in the
context the Torelli group.  Bestvina-Bux-Margalit exploited the action of the Torelli group on this complex to
compute the dimension of the Torelli group.  Hatcher-Margalit used it to identify generating sets for the Torelli group.

In \cite{II}, Irmer defines the homology curve complex of a surface:

``Suppose $S$ is a closed oriented surface.  $S$ is not required to be connected but every component is
assumed to have genus $g \geq 2$.

Let $\alpha$ be a nontrivial element of $H_1(S, {\mathbb Z})$.  The {\em homology curve complex},
${\cal HC}(S, \alpha)$, is a simplicial complex whose vertex set is the set of all homotopy classes of
oriented multi-curves in $S$ in the homology class $\alpha$.  A set of vertices $m_1, \dots, m_k$ spans
a simplex if there is a set of pairwise disjoint representatives of the homotopy classes.

The {\em distance}, $d_{\cal H}(v_1, v_2)$, between two vertices $v_1$ and $v_2$ is defined to
be the distance in the path metric of the one-skeleton, where all edges have length one."  
(See \cite[Page 1]{II}.)

It is not hard to see the following ({\it cf,} Remark \ref{simplex} and Figure \ref{ccstar}):

\begin{lem}
When both are defined, {\it i.e.,} when $S$ is closed, connected, of genus at least $2$ and
$\alpha$ is primitive, $Kak(S, \alpha)$ is a subcomplex of ${\cal HC}(S, \alpha)$.   Moreover, for vertices
$v, v'$ of $Kak(S, \alpha)$,
\[d_K(v, v') \geq d_{\cal H}(v, v')\]
\end{lem}

Irmer shows that distance between vertices of ${\cal HC}(S, \alpha)$ is bounded above by a linear function 
on the intersection number of representatives.   The same is true for vertices of the Kakimizu complex.   Irmer
also constructs quasi-flats in ${\cal HC}(S, \alpha)$.  Her construction carries over to the setting of the 
Kakimizu complex.  See Section \ref{quasi-flats}.

\section{The projection map, distances and geodesics} \label{pdg}

In \cite{K}, Kakimizu defined a map on the vertices of the Kakimizu complex of a knot.
He used this map to prove several things, for instance that the metric, $d_K$, on the
vertices of the Kakimizu complex equals graph distance.  (Quoted and reproved here as
Theorem \ref{distance}.)
In \cite{PS}, Kakimizu's map was rebranded as a projection map.  









We wish to define 
\[\pi_{Vert(Kak(S, \alpha))}:Vert(Kak(S, \alpha)) \rightarrow Vert(Kak(S, \alpha))\]
on the vertex set of $Kak(S, \alpha)$.  
Let $(p, \hat S, S)$ be the infinite cyclic covering space associated with $\alpha$.  
Let $v, v'$ be vertices in $Kak(S, \alpha)$ such that $v \neq v'$.  
Here $v = [(w, c)]$ for some compact oriented 1-manifold $c$ and $v' = [(w', c')]$ 
for some compact oriented 1-manifold $c'$.  We may assume, in accordance with Definition
\ref{d_K} and Remark 
\ref{d_Kisfinite}, that $(w, c)$ and $(w', c')$ are chosen so that
$d_K(c, c') = d_K(v, v')$.  Define $\tau, S_i, S'_i, c_i$ and, by analogy, $c'_i$, as in 
Definition \ref{d_K}.  

\begin{figure}[p]
\vspace{2 mm}
\centerline{\epsfxsize=2.5in \epsfbox{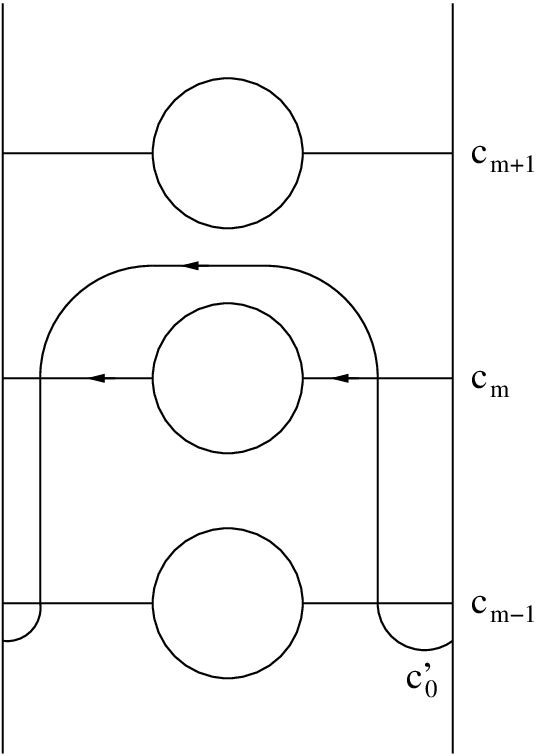}}
\caption{\sl The setup with $c_m, c'_0$ (weights are 1)}
\label{proj}
\vspace{2 mm}
\end{figure}

Instead of working only with $c'_0$, we will now also work with $h(w', c'_0)$.   
Take $m = max \{ i \; | \; S_{i+1} \cap S'_0 \neq \emptyset \}$.  
Consider a connected component $C$ of $S_{m+1} \cap S'_0$.  Its frontier consists of
a subset of $c'_0$ and a subset of $c_m$.  The subset of $c'_0$ lies above the subset
of $c_m$.  In particular, $C$ lies above $c_m$ and below $c'_0$, hence the 
orientations of the subset of $c'_0$ are opposite those of the
subset of $c_m$.  See Figure \ref{proj}.  Because the subset of $c'_0$ and the
subset of $c_m$ cobound $C$, they are homologous.  It follows that the lowest components of the corresponding subset of 
$h(w', c'_0)$ are also homologous to the subset of $c_m$.  

Replacing the lowest of the corresponding subsets of $h(w', c'_0)$ with the subset of $c_m$ and isotoping this portion of 
$c_m$ to lie below $c_m$
yields a multi-curve $d_1$ with the following properties: 

\begin{itemize}

\item
$d_1$ is homologous to $h(w', c_0')$ via a homology that descends to a homology in $S$ (because $C$ is homeomorphic to a subset of $S$); 

\item
$d_1$ has lower geometric intersection number with $c_m$ than $h(w', c_0')$; 

\item
$d_1$ lies above $h(w', c_{-1}')$ and can 
be isotoped to lie below and thus be disjoint from $h(w', c_0')$, moreover its projection can be isotoped to be disjoint from $h(w', c')$.  

\item
For $(x_1,e_1) = h^{-1}(d_1)$, we have $x_1 \circ e_1$ homologous to $w' \circ c_0'$ via
a homology that descends to a homology in $S$.

\end{itemize}

\begin{figure}[p]
\vspace{2 mm}
\centerline{\epsfxsize=2.5in \epsfbox{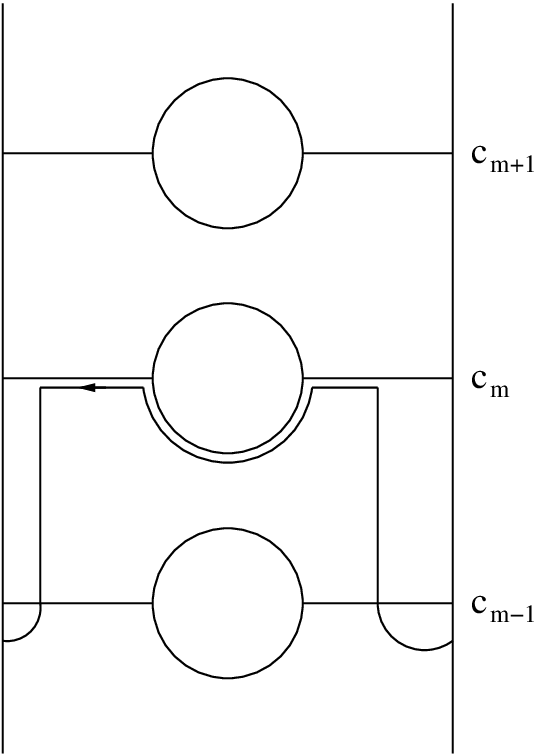}}
\caption{\sl $d_1$}
\label{afterproj}
\vspace{2 mm}
\end{figure}

See Figures \ref{proj}, \ref{afterproj}.
 
\begin{figure}[p]
\vspace{2 mm}
\centerline{\epsfxsize=2.5in \epsfbox{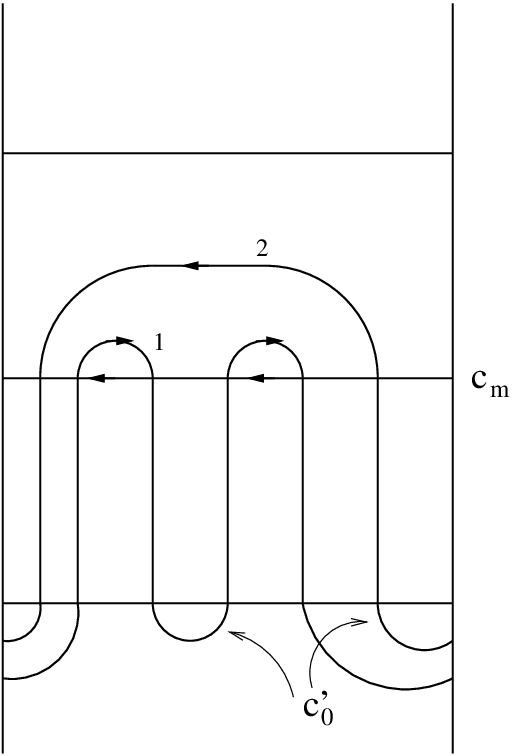}}
\caption{\sl A different pair of weighted multi-curves}
\label{c0cm}
\vspace{2 mm}
\end{figure}

\begin{figure}[p]
\vspace{2 mm}
\centerline{\epsfxsize=2.5in \epsfbox{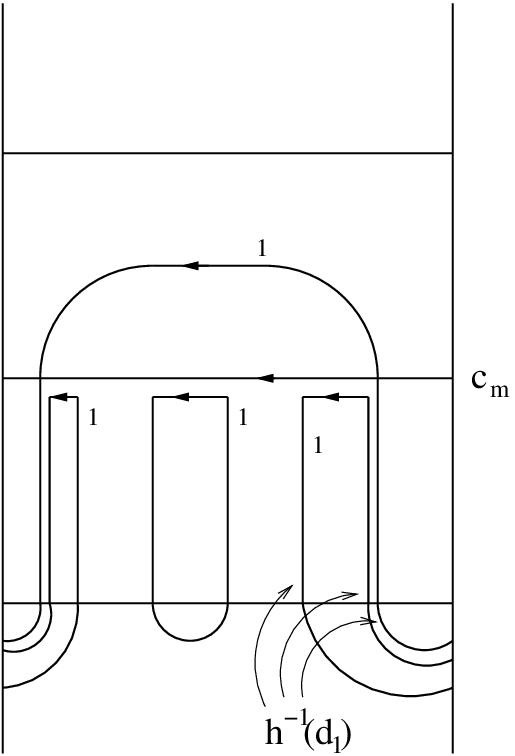}}
\caption{\sl $h^{-1}(d_1)$}
\label{d1}
\vspace{2 mm}
\end{figure}

Working with $h^{-1}(d_1), d_1$ instead of $c_0', h(w', c_0')$ we perform such replacements in succession to obtain a sequence of 
multi-curves $d_1, \dots, d_k$ such that the following hold:

\begin{itemize}

\item
$d_j$ is homologous to $d_{j-1}$ via a homology that descends to a homology in $S$; 

\item
$d_j$ has lower geometric intersection number with $c_m$ than $d_{j-1}$; 

\item
$d_j$ can 
be isotoped to lie below $h(w', c_0'), d_1, \dots, d_{j-1}$, moreover its projection can be isotoped to be disjoint from $h(w', c')$.  

\item
For $(x_j,e_j) = h^{-1}(d_j)$, we have $x_j \circ e_j$ homolgous to $w' \circ c_0'$ via
a homology that descends to $S$.

\item
$d_k$ lies above $h(w', c_{-1}')$ and below $c_m$.

\end{itemize}

See Figures \ref{c0cm}, \ref{d1} and \ref{d2}.  This proves the following:

\begin{figure}[h]
\vspace{2 mm}
\centerline{\epsfxsize=2.5in \epsfbox{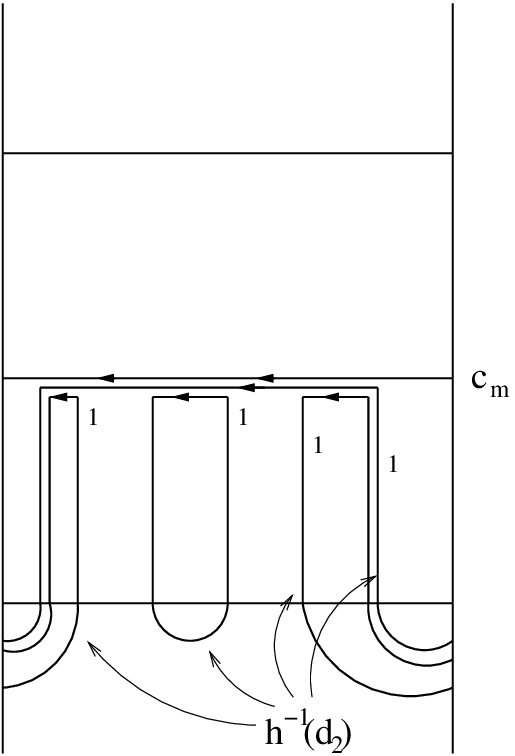}}
\caption{\sl $h^{-1}(d_2)$}
\label{d2}
\vspace{2 mm}
\end{figure}

\begin{lem} \label{pucc}
The homology class $[[p(d_k)]] = p_{\#}(x_k \circ e_k) = \alpha$.
\end{lem}

We make two observations: 1) A result of Oertel, see \cite{Oe1988}, shows that the isotopy class of 
$p(e_k)$ does not depend on the choices 
made; 2) It is important to realize that $(x_k, p(e_k))$ may not be a Seifert curve, because
$S \backslash p(e_k)$ is not necessarily connected.

\begin{figure}[h]
\vspace{2 mm}
\centering
\includegraphics[width=2.5in]{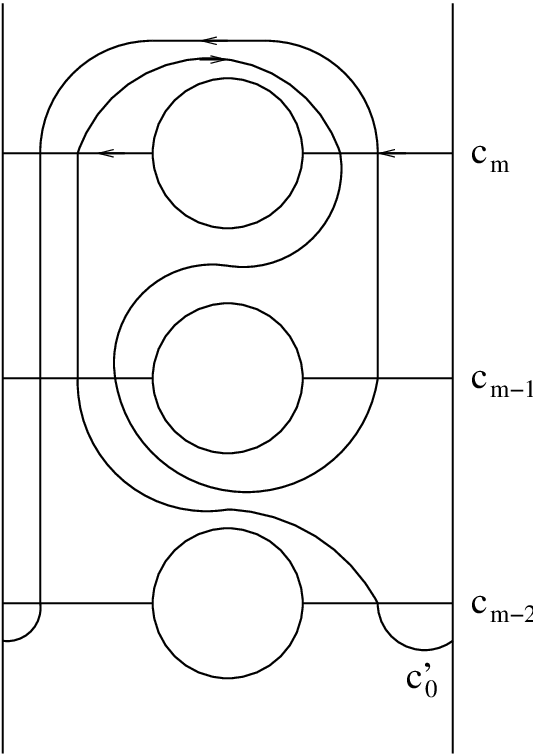}
\caption{\sl The setup with $c_m, c_0'$ (weights are 1)}
\label{projsnail}
\vspace{2 mm}
\end{figure}
\begin{figure}[p]
\vspace{2 mm}
\centerline{\epsfxsize=2.5in \epsfbox{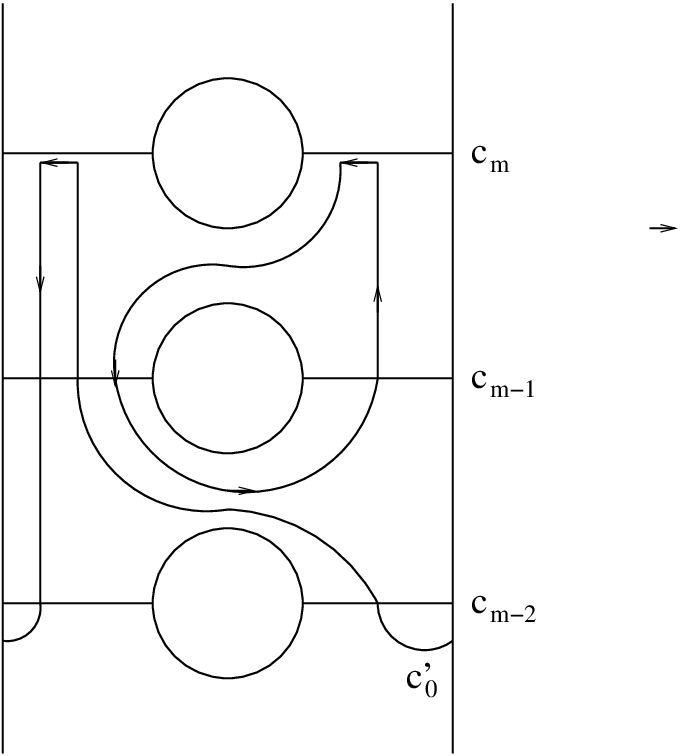}}
\caption{\sl $d_k$}
\label{afterprojsnail}
\vspace{2 mm}
\end{figure}
\begin{figure}[p]
\vspace{2 mm}
\centerline{\epsfxsize=2.5in \epsfbox{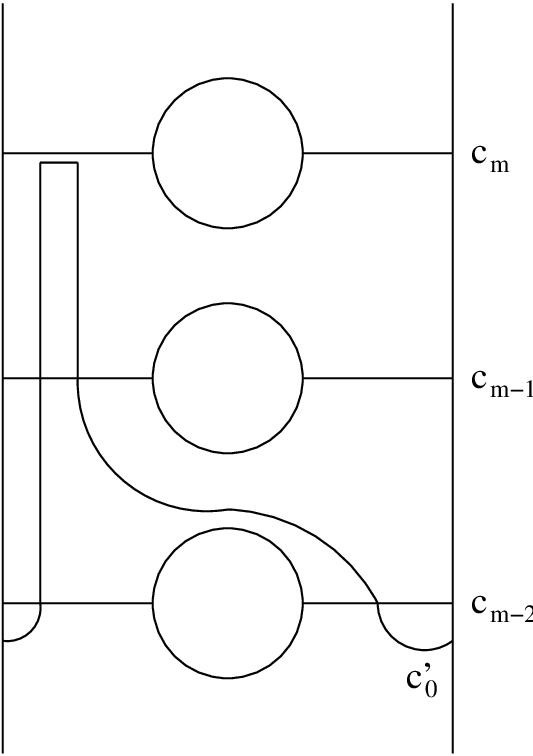}}
\caption{\sl A subset of $d_k$}
\label{afterprojsnailafter}
\vspace{2 mm}
\end{figure}

If $S \backslash p(h^{-1}(e_k))$
is connected, set $p_c(c') = (x_k, p(e_k))$.  
Otherwise, choose a component $D$ of 
$S \backslash p(e_k)$.  If the frontier of $D$ is null homologous, then remove the frontier of $D$ from $p(e_k)$.
See Figures \ref{projsnail}, \ref{afterprojsnail}, \ref{afterprojsnailafter}.

If the frontier of $D$ is not null homologous (because the orientations do not match up) choose an arc $a$ in its frontier with smallest weight.
Denote the weight of $a$ by $w^a$.  We eliminate the component $a$ of $p(e_k)$ by adding $\pm w^a$ to
the weights of the other components of  $p(e_k)$ in the frontier of $D$ in such a way that
the resulting weighted multi-curve still has homology $\alpha$.      

After a finite number of such eliminations, we obtain a weighted multi-curve that is a subset of 
$p(e_k)$, has homology $\alpha$, and whose complement in $S$ is connected.  After reversing orientation on components
with negative weights, we obtain a Seifert curve $p_c(c')$.  

\begin{lem} \label{ucc}
The homology class $[[p_c(c')]] = \alpha$.
\end{lem}

\proof
This follows from Lemma \ref{pucc} and the observations above.
\qed


\begin{defn}
We denote the isotopy class $[p_c(c')]$ by $\pi_v(v')$.
\end{defn}

\begin{lem} \label{ineqs}
For $v \neq v'$, the following hold:
\[d_K(\pi_v(v'), v') = 1\] and \[d_K(\pi_v(v'), v) \leq d_K(v', v) - 1.\]
\end{lem}

It will follow from Theorem \ref{distance} below that the inequality is in fact an equality.  

\proof
By construction, $e_k$ lies strictly between $c_0'$ and $c_{-1}'$.  
So $\tau(e_k)$ lies strictly between $c_1'$ and $c_0'$.  Thus
the lift of $S \backslash p_c(c')$ with frontier in
$e_k \cup \tau(e_k)$
meets $S_0'$ and $S_1'$ and is disjoint from $S_i'$ for $i \neq 0, 1$.  
It follows that the lift of $S \backslash p_c(c')$ with frontier contained in $e_k \cup \tau(e_k)$ also meets 
$S_0'$ and $S_1'$ and is disjoint from $S_i'$ for 
$i \neq 0, 1$.   
Hence \[d_K(\pi_v(v'), v') = 1.\]

In addition, suppose that $c_0' \cap S_i \neq \emptyset$ if and only if $i \in \{n, \dots, m+1\}$. 
Then $c_1' \cap S_i \neq \emptyset$ if and only if
$i \in \{n+1, \dots, m+2\}$.  
Hence the lift of 
$S \backslash c'$ that lies strictly between $c_0'$ and $c_1'$ meets exactly 
$S_n, \dots, S_{m+2}$.  

By construction, 
$e_k \cap S_i$ can be non empty only if $i \in \{n, \dots, m\}$ and
thus $\tau(e_k) \cap S_i$ can be non empty only if $i \in \{n+1, \dots, m+1\}$.
Hence the lift of 
$S \backslash p_c(c')$ with frontier in $e_k \cup \tau(e_k)$ 
can meet $S_i$ only if $i \in \{n, \dots, m+1\}$.  It follows that the lift of $S \backslash p_c(c')$ with frontier contained in 
$e_k \cup \tau(e_k)$ 
can meet $S_i$ only if $i \in \{n, \dots, m+1\}$.
Whence \[d_K(\pi_v(v'), v) \leq m + 1 - n - 1  = d_K(v', v) - 1.\]
\qed

\begin{defn}
The {\em graph distance} on a complex ${\cal C}$ is a function that assigns 
to each pair of vertices $v, v'$ the least possible number of edges in an 
edge path in ${\cal C}$ from $v$ to $v'$.  
\end{defn}

\begin{thm}(Kakimizu) \label{distance}
The function $d_K$ equals graph distance.
\end{thm}

\proof
Denote the graph distance between $v'$ and $v$ by $d(v', v)$.
If $d_K(v', v) = 1$, then $d(v', v) = 1$ and vice versa by definition.
So suppose $d_K(v', v) = m > 1$ and consider the path with vertices \[v', \pi_v(v'), \pi^2_v(v'), \dots, 
\pi_v^{m-1}(v'), \pi_v^m(v') = v.\]
By Remark \ref{ineqs}, $d_K(\pi_v(v'), v') = 1$ and $d_K(\pi_v^i(v'), \pi_v^{i-1}(v')) = 1$.
Thus the existence of this path guarantees that $d(v', v) \leq m$.  
Hence $d(v',v) \leq d_K(v',v)$.  
Let $v' = v_0, v_1, \dots, v_n = v$ be the vertices of 
a path realizing $d(v', v)$.
By the triangle inequality and the fact that $d(v_{i-1}, v_i) = 1 = d_K(v_{i-1}, v_i)$:

\[d_K(v', v) \leq d_K(v_0, v_1) + \dots + d_K(v_{n-1}, v_n) = 1 + \dots + 1 =\]
\[d(v_0, v_1) + \dots + d(v_{n-1}, v_n) = d(v', v)\]
\qed

The following theorem is a reinterpretation of a theorem of Scharlemann and Thompson, see \cite{ST}, that was proved using different methods:

\begin{thm}
The Kakimizu complex is connected.
\end{thm}

\proof
Let $v, v'$ be vertices in $Kak(S, \alpha)$.  By Remark \ref{d_Kisfinite}, $d_K(v, v')$ is finite.  
By Theorem \ref{distance}, $d(v, v')$ is
finite.  In particular, there is a path between $v$ and $v'$.
\qed

\begin{defn}
A {\em geodesic} between vertices $v, v'$ in a Kakimizu complex is an edge-path that realizes
$d(v, v')$.  
\end{defn}

\begin{thm} \label{geodesics}
The path with vertices $v', \pi_v(v') , \pi_v^2(v'), \dots, \pi_v(v')^{m-1}, \pi^m_v(v') = v$ is a geodesic.
\end{thm}

\proof
This follows from Theorem \ref{distance} because the path
\[v', \pi_v(v') , \pi_v^2(v'), \dots, \pi_v(v')^{m-1}, \pi^m_v(v') = v\]
realizes $d(v', v)$.
\qed

\begin{rem}
Theorem \ref{geodesics} tells us that geodesics in the Kakimizu complex joining two given vertices 
are, at least theoretically, constructible.  
\end{rem}

\begin{figure}[ht]
\vspace{2 mm}
\centerline{\epsfxsize=2.5in \epsfbox{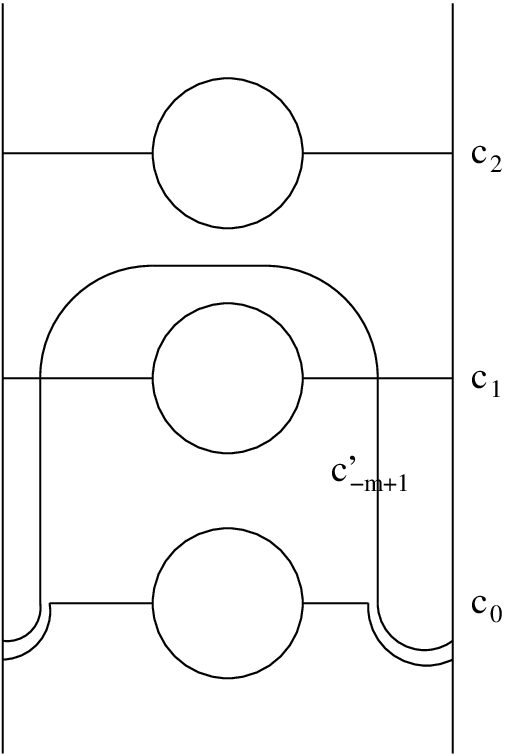}}
\caption{\sl $u(c_{-m+1}', c_0)$}
\label{afterprojalt}
\vspace{2 mm}
\end{figure}

Note that, typically, $\pi_v(v') \neq \pi_{v'}(v)$.  See Figure \ref{afterprojalt} for a step in the
construction of $\pi_{v'}(v)$.  

\section{Contractibility} \label{contractibility}

The proof of contractibility presented here is a streamlined version of the proof given in the 
3-dimensional case in \cite{PS}.  Those familiar with Hatcher's work in \cite{Hatcher}, will note 
certain similarities with his first proof of contractibility of $C_{\alpha}(S)$ in the case 
that $\alpha$ is primitive.  

\begin{lem} \label{flag}
Suppose that $v, v^1, v^2$ are vertices in $Kak(S, \alpha)$.  Then there are representatives
$c, c^1, c^2$ with $v = [(w, c)],$  $v^1 = [(w^1, c^1)],$ and $v^2 = [(w^2, c^2)]$ that realize 
$d_K(v, v^1), d_K(v, v^2),$ and $d_K(v^1, v^2)$.
\end{lem}

\proof
Let $c, c^1, c^2$ be geodesic representatives of the underlying curves of representatives of 
$v, v^1, v^2$ such that 
arc components of $c, c^1, c^2$ are perpendicular to $\partial S$.
Lifts of $c, c^1,$ and 
$c^2$ to $(p, \hat S, S)$, the infinite cyclic covering of $S$ associated with $\alpha$, are also 
geodesics.  Points of intersection lift to points of intersection.  
Geodesics that intersect can't be isotoped to be disjoint.  Hence $c, c^1, c^2$, with appropriate weights, 
realize $d_K(v, v^1),$ $d_K(v, v^2),$ and $d_K(v^1, v^2)$. 
\qed

\begin{lem} \label{edges}
Suppose that $v, v^1, v^2$ are vertices in $Kak(S, \alpha)$ such that $d_K(v, v^i) > 1$ and $d_K(v^1, v^2) = 1$.  Then
$d_K(\pi_v(v^1), \pi_v(v^2)) \leq 1$.  
\end{lem}

\proof
In the case that, say, $v^1 = v$, note that $d_K(v^1, v^2) = 1$ means that $d_K(v, v^2) = 1$.  
Thus $\pi_v(v^1) = \pi_v(v) = v$ and $\pi_v(v^2) = v$.   Thus $d_K(\pi_v(v^1), \pi_v(v^2)) = 0$.
In the case that, say, $d_K(v^1, v) = 1$ and $v^2 \neq v$, note that $\pi_v(v^1) = v$ and 
\[d_K(v, v^2) \leq d_K(v, v^1) + d_K(v^1, v^2) = 1 + 1\]
thus \[d_K(v, \pi_v(v^2) \leq 1,\] by Lemma \ref{ineqs}, and  $d_K(\pi_v(v^1), \pi_v(v^2)) \leq 1$.
Hence we will assume, for the rest of this proof, that $d_K(v, v^i) > 1$.

By Lemma \ref{flag}, there are
representatives $(w, c), (w^1, c^1),$ and $(w^2, c^2)$ of $v, v^1$ and $v^2$ that realize
$d_K(v, v^1), d_K(v, v^2),$ and $d_K(v^1, v^2)$. 
Let $(p, \hat S, S)$ be the infinite cyclic cover of $S$ 
associated with $\alpha$.  Define $\tau, S_i, S_i^1, S_i^2, c_i, c_i^1, c_i^2$ as in 
Definition \ref{d_K} but with a caveat: Label $S_i^1, S_i^2$ so that $S_0^1, S_0^2$
meet $S_1$ and meet $S_j$ only if $j \leq 1$.

\begin{figure}[p]
\vspace{2 mm}
\centerline{\epsfxsize=2.5in \epsfbox{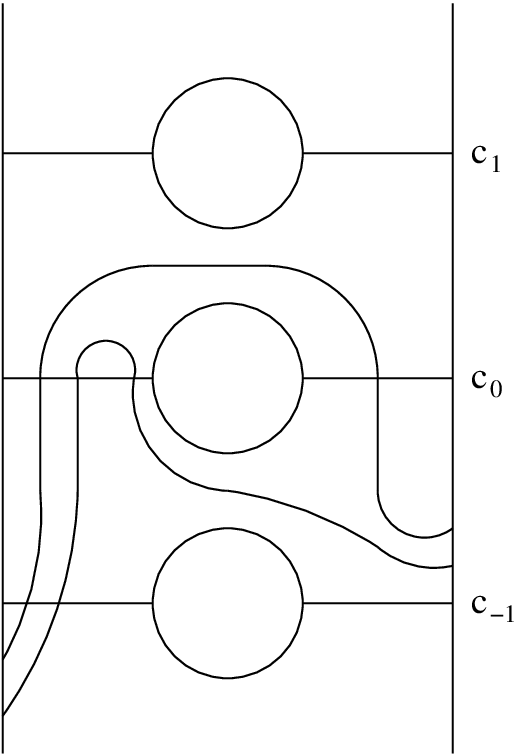}}
\caption{\sl $c_0^1$ and $c_0^2$}
\label{2xproj}
\vspace{2 mm}
\end{figure}

\begin{figure}[p]
\vspace{2 mm}
\centerline{\epsfxsize=2.5in \epsfbox{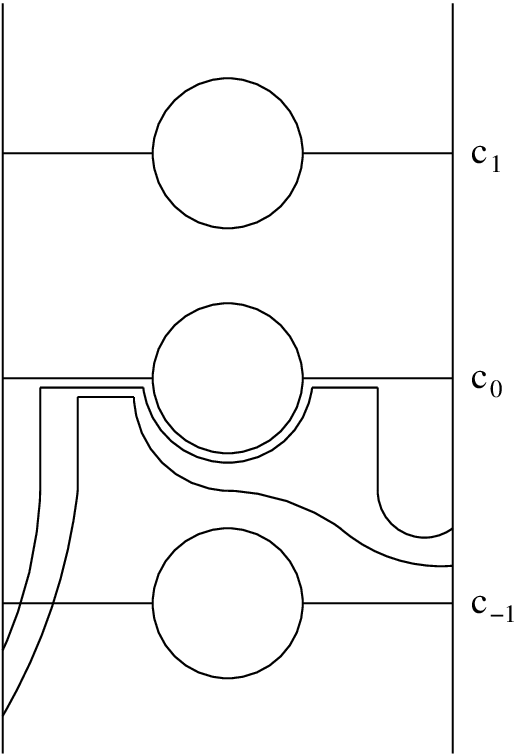}}
\caption{\sl $d_k^1$ and $d_l^2$}
\label{2xafterproj}
\vspace{2 mm}
\end{figure}

Since $d_K(c^1, c^2) = d_K(v^1, v^2)$, $c^1$ and $c^2$ must be disjoint.  
Since $c_0^1$ is separating,
$c_0^2$ lies either above or below $c_0^1$.  Without loss of generality, we will assume
that $c_0^2$ lies above $c_0^1$ (and below $\tau(c_0^1)$).   See Figures \ref{2xproj} and \ref{2xafterproj}.
Note that $h(w^2, c_0^2)$ also lies above $h(w^1, c_0^1)$.  Proceeding as in the discussion preceding Lemma \ref{ucc}, 
construct $e_k^1$ whose projection contains $p_c(c^1)$ and then $e_l^2$ whose projection contains $p_c(c^2)$,
noting that this
construction can be undertaken so that $e_l^2$ lies above 
$e_k^1$ (and below $\tau(e_k^1)$).  

Consider the lift of $S \backslash p_c(c^2)$
with frontier in $e_l^2 \cup \tau(e_l^2)$.  This lift
of $S \backslash p_c(c^2)$
meets at most the two lifts of $S \backslash p_c(c^1)$ whose frontiers lie in 
$e_k^1 \cup \tau(e_k^1)$ and 
$\tau(e_k^1) \cup \tau^2(e_k^1)$.
Whence 
\[d_K(\pi_v(v^1), \pi_v(v^2)) \leq 1.\]
\qed



\begin{lem} \label{continuous}
If $d_K(v^1, v^2) = m$, then $d_K(\pi_v(v^1), \pi_v(v^2)) \leq m$.
\end{lem}

\proof
Let $v^1 = v_0, v_1, \dots, v_{m-1}, v_m = v^2$ be the vertices of a path from $v^1$ to $v^2$ 
that realizes $d_K(v^1, v^2)$.  By Lemma \ref{edges}, 
$d_K(\pi_v(v_i), \pi_v(v_{i+1})) \leq 
d_K(v_i, v_{i+1}) = 1$ for $i = 0, \dots, m-1$.
Hence \[d_K(\pi_v(v^1), \pi_v(v^2)) \leq d_K(\pi_v(v_0), \pi_v(v_1)) + \dots + d_K(\pi_v(v_{m-1}), \pi_v(v_m)) \leq\]
\[d_K(v_0, v_1) + \dots + d_K(v_{m-1}, v_m) \leq m.\]
\qed

\begin{thm}
The Kakimizu complex of a surface is contractible.
\end{thm}

\proof
Let $Kak(S, \alpha)$ be a Kakimizu complex of a surface.  It is well known (see \cite[Exercise 11, page 358]{HatcherAT}) that
it suffices to show that every finite subcomplex of $Kak(S, \alpha)$ is contained in a contractible 
subcomplex of $Kak(S, \alpha)$.  Let ${\cal C}$ be a finite subcomplex of $Kak(S, \alpha)$.
Choose a vertex $v$ in ${\cal C}$ and denote by ${\cal C}'$ the smallest flag complex
containing every geodesic of the form given in Theorem \ref{geodesics} for $v'$ a vertex 
in ${\cal C}$.  Since ${\cal C}$ is finite, it follows that ${\cal C}'$ is finite.  

Define 
$c:Vert({\cal C}') \rightarrow Vert({\cal C}')$ on vertices by
$c(v') = \pi_v(v')$.  By Lemma \ref{edges}, this map extends to edges.  Since 
${\cal C}'$ is flag, the map extends to simplices and thus to all of
${\cal C}'$.  By Lemma \ref{continuous} this map is continuous.  It is not hard to see that 
$c$ is homotopic to the identity map.  In particular, $c$ is a contraction map.  
(Specifically, $c^d$, where $d$ is the diameter of ${\cal C}'$, has the set $\{v\}$ as its image.)
\qed

\section{Dimension}

In \cite{Hatcher}, Hatcher proves that the dimension of $C_{\alpha}(S)$ is $2g(S) - 3$, where
$g(S)$ is the genus of the closed oriented surface $S$.  An analogous argument derives the
same result in the context of $Kak(S, {\alpha})$.    

\begin{lem} \label{closedkakdim}
Let $S$ be a closed connected orientable surface with $genus(S) \geq 2$ and $\alpha$ a primitive class in 
$H_1(S, \partial S)$.  The dimension of $Kak(S, \alpha)$ is $-\chi(S) - 1 = 2genus(S) - 3$.
\end{lem}

\begin{figure}[h]
\vspace{2 mm}
\centerline{\epsfxsize=5in \epsfbox{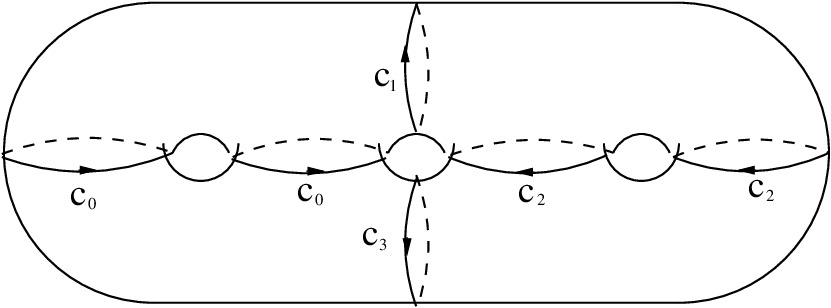}}
\caption{\sl A simplex in a genus $3$ surface}
\label{2g-3}
\vspace{2 mm}
\end{figure}

\proof
It is not hard to build a simplex of $Kak(S, \alpha)$ of dimension $2genus(S) - 3$.
See for example Figure \ref{2g-3}, where $0, 1, 2,$ and  $3$ are multi-curves (each of weight $1$) 
representing the vertices of a simplex.   Thus the dimension of $Kak(S, \alpha)$ is
greater than or equal to $2genus(S)- 3$.

Conversely, let $\sigma$ be a simplex of maximal dimension in $Kak(S, \alpha)$.  Label the vertices
of $\sigma$ by $v_0, \dots, v_n$ and let $c_0, \dots, c_n$ be geodesic representatives of the underlying curves 
of representatives of 
$v_0, \dots, v_n$.
By Lemma \ref{cycleofcycles},
 $S \backslash (c_0 \cup \dots \cup c_n)$ consists of subsurfaces $P_0, \dots, P_n$ with frontiers
$c_0 - c_n, c_1 - c_0, \dots, c_n - c_{n-1}$.  Since $c_i$ and $c_{i-1}$ are not isotopic, no $P_i$ can
 consist of annuli.  In addition, no $P_i$ can be a sphere, hence each must have negative Euler characteristic.
Thus the number of $P_i$'s is at most $-\chi(S)$.  
{\it I.e.,} \[n \leq - \chi(S) = 2genus(S) -2.\]  In other words, the dimension of $\sigma$ and hence the
dimension of $Kak(S, \alpha)$ is less than or equal to $2genus(S) - 3$.
\qed

We can extend this argument to compact surfaces, by introducing the following notion of complexity:

\begin{defn}
Let $S$ be a compact surface and let $P$ be an open subset of $S$ whose boundary
consists of open subarcs of $\partial S$ and, possibly, 
components of $\partial S$.  Define
\[c(P, S) = -2\chi(P) + number\;of\;open\;subarcs\;in\;\partial P\]
\end{defn}

The following lemma is immediate:

\begin{lem}
Let ${\cal C}$ be a union of simple closed curves and simple arcs in $S$.  Then
\[c(S, S) = c(S \backslash {\cal C}, S)\]
\end{lem}

\begin{thm} \label{kakdim}
Let $S$ be a compact connected orientable surface with $\chi(S) \leq -1$ and $\alpha$ a primitive class in 
$H_1(S, \partial S)$.  The dimension of $Kak(S, \alpha)$ is $-2\chi(S) - 1 = 4genus(S) + 2b - 5$, where
$b$ is the number of boundary components of $S$.
\end{thm}

\begin{figure}[h]
\vspace{2 mm}
\centerline{\epsfxsize=4in \epsfbox{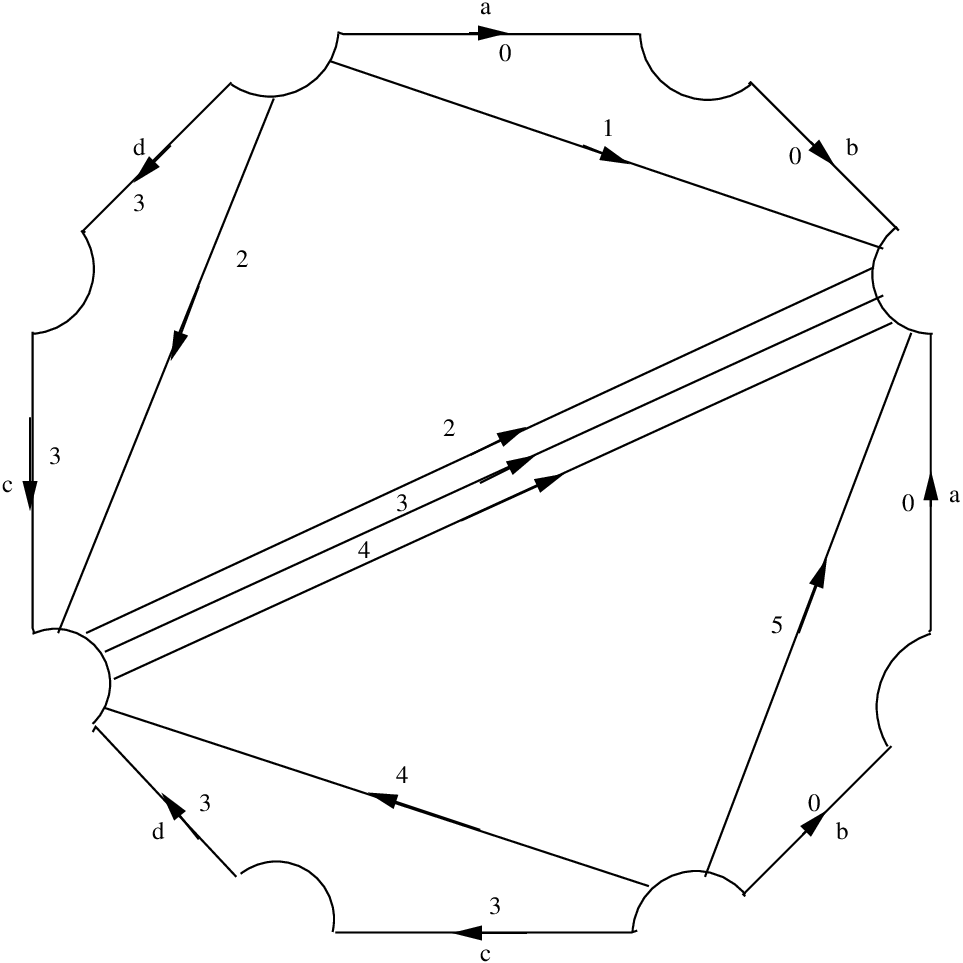}}
\caption{\sl A simplex in a punctured genus $2$ surface}
\label{4g+2b-5}
\vspace{2 mm}
\end{figure}

\proof
To build a simplex of $Kak(S, \alpha)$ of dimension $4genus(S) + 2b - 5$, 
see for example Figure \ref{4g+2b-5}, where $0, 1, 2, 3, 4$ and  $5$ are multi-curves (each of weight $1$) 
representing the vertices of a simplex.    Thus the dimension of $Kak(S, \alpha)$ is
greater than or equal to $4g(S)+2b-5$.

Conversely, let $\sigma$ be a simplex of maximal dimension in $Kak(S, \alpha)$.  Label the vertices
of $\sigma$ by $v_0, \dots, v_n$ and let $c_0, \dots, c_n$ be geodesic representatives of the underlying curves of representatives of $v_0, \dots, v_n$ such that 
arc components of $c_0, \dots, c_n$ are perpendicular to $\partial S$.
By Lemma \ref{cycleofcycles}, $S \backslash (c_0 \cup \dots \cup c_n)$ consists of subsurfaces 
$P_0, \dots, P_n$ with frontiers containing
$c_0 - c_n, c_1 - c_0, \dots, c_n - c_{n-1}$.
Since $c_i$ and $c_{i-1}$ are not isotopic, $P_i$ can't
 consist of annuli or disks with exactly two open subarcs of $\partial S$ in their boundary.  In addition, no $P_i$ can be 
a sphere or a disk with exactly
one open subarc of $\partial S$ in its boundary, hence each must have positive complexity.
Thus the number of $P_i$'s is at most $c(S,S\backslash (c_0 \cup \dots \cup c_n))$.  
{\it I.e.,} \[n \leq c(S, S \backslash (c_0 \cup \dots \cup c_n)) = c(S,S) = - 2\chi(S).\]  In other words, 
the dimension of $\sigma$ and hence the
dimension of $Kak(S, \alpha)$ is less than or equal to $-2\chi(S) - 1 = 4genus(S) + 2b - 5$.
\qed

\section{Quasi-flats} \label{quasi-flats}

In this section we explore an idea of Irmer.  See \cite[Section 7]{II}.  

\begin{figure}[ht]
\vspace{2 mm}
\centerline{\epsfxsize=5in\epsfbox{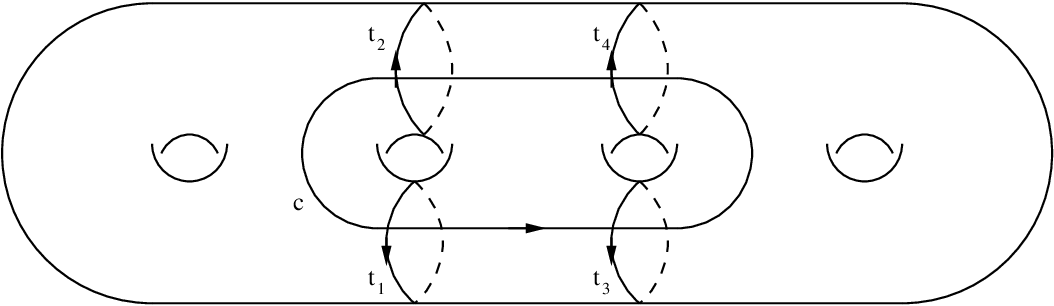}}
\caption{\sl Building a quasi-flat by Dehn twists}
\label{irmerexample}
\vspace{2 mm}
\end{figure}

Consider Figure \ref{irmerexample}.  Denote the surface depicted by $S$ and the homology class of $c$ by $\alpha$. 
The curves $t_1$ and $t_2$ are homologous as are $t_3$ and $t_4$.
Denote by $v$ the vertex $(1, c)$  of $Kak(S, \alpha)$, 
by $v_1$ the vertex corresponding to the result, $d_1$, obtained from $c$ by Dehn twisting $n$  
times around $t_1$ and $-n$ times around $t_2$, and by $v_2$ the vertex corresponding to the result, $d_2$,
obtained from $c$ by Dehn twisting $n$ times around $t_3$ and $-n$ times around $t_4$.   
Then $d_1, d_2$ are homologous to $c$, so we obtain three vertices $v, v_1, v_2$ in $Kak(S, \alpha)$ (all weights are 1).   
Note the following:

\[d(v, v_i) = d_K(v, v_i) = n\] 
\[d(v_1, v_2) = d_K(v_1, v_2) = n\]

For $i = 1, 2,$ we consider the geodesics $g_i$ with vertices $v_i, \pi_v(v_i), \dots, \pi_v^n(v_i) = v$.  
In addition, consider the geodesic $g_3$ with vertices $v_2, \pi_{v_1}(v_2), \dots, \pi_{v_1}^n(v_2) = v_1$ and
note that $\pi_{v_1}^i(v_2)$ is represented by a curve obtained from $c$ by Dehn twisting $i$ times around 
$t_1$, $-i$ times around $t_2$,  $n-i$  times around $t_3$ and $-(n-i)$ times around $t_4$.

\begin{defn} 
Let $(X, d)$ be a metric space.  A {\em triangle} is a 6-tuple $(v^1,$ $v^2,$ $v^3,$ $g^1,$ $g^2,$ $g^3)$, where
$v^1, v^2, v^3$ are vertices and the edges $g^1, g^2, g^3$ satisfy the following: $g^1$ is a distance 
minimizing path between $v^1$ and $v^2$, $g^2$ is a distance minimizing path between $v^2$ and $v^3$,
$g^3$ is a distance minimizing path between $v^3$ and $v^1$.  

A triangle $(v^1, v^2, v^3, g^1, g^2, g^3)$ 
is {\em $\delta$-thin} if each $g^i$ lies in a $\delta$-neighborhood of the other two edges.
A metric space $(X, d)$ is {\em $\delta$-hyperbolic} if every triangle in $(X, d)$ is $\delta$-thin. 
It is {\em hyperbolic} if there is a $\delta > 0$ such that $(X, d)$ is $\delta$-hyperbolic. 
\end{defn}

For $n$ even, the midpoint, $m_1$, of the geodesic $g_1$ is the 
vertex corresponding to the result, $d_1'$, obtained from $c$ by Dehn twisting $\frac{n}{2}$ times around 
$t_1$ and $-\frac{n}{2}$ times around $t_2$.  Likewise, the midpoint, $m_2$, of the geodesic $g_2$ 
is the vertex corresponding to the result, $d_2'$, obtained from $c$ by Dehn twisting 
$\frac{n}{2}$ times around $t_3$ and $-\frac{n}{2}$ times around $t_4$.  
The midpoint, $m_3$, of $g_3$ is represented by a curve obtained from $c$ by Dehn twisting 
$\frac{n}{2}$ times around $t_1$ and around $t_3$ and $-\frac{n}{2}$ times around $t_2$ and $t_4$.  

\begin{lem} \label{nonhyp}
Let $S$ be the closed oriented surface of genus $4$.  Then $Kak(S, \alpha)$ is not 
hyperbolic.  
\end{lem}

\proof
For $S$ the closed genus $4$ surface, the triangle $(v, v_1, v_2, g_1, g_2, g_3)$ described depends on $n$, 
so we will denote it by $T_n$.
In $T_n$ we have the following:

\[d(v, m_3) = d_K(v, m_3) = n\] 
\[d(v_1, m_2) = d_K(v_1, m_2) = n\]
\[d(v_2, m_1) = d_K(v_2, m_1) = n\]

In particular, $g_3$ is contained in a $\delta$-neighborhood of the two geodesics $g_1, g_2$
only if $n$ is less than $\delta$.   
Thus the triangle $T_n$ in $Kak(S, \alpha)$ is not
$\delta$-thin for $n \geq \delta$.  It follows that $Kak(S, \alpha)$ is not hyperbolic.  
\qed

\begin{defn}
Let $(X, d)$ be a metric space.  A {\em quasi-flat} in $(X, d)$ is a quasi-isometry from ${\mathbb R}^n$ to $(X, d)$,
for $n \geq 2$.   
\end{defn}

Note the following:

\[d(m_1, m_2) = d_K(m_1, m_2) = \frac{n}{2}\] 
\[d(m_1, m_3) = d_K(m_1, m_3) = \frac{n}{2}\]
\[d(m_2, m_3) = d_K(m_2, m_3) = \frac{n}{2}\]

Thus the triangle $T_n$ scales like a Euclidean triangle.  It is not too hard to see that a triangle with this property
can be used to construct a quasi-isometry between ${\mathbb R}^2$ and an infinite union of such triangles lying in
$Kak(S, \alpha)$.  Thus $Kak(S, \alpha)$ contains quasi-flats.  It is also not hard to adapt this construction to show
that, for $S$ an oriented surface, $Kak(S, \alpha)$ is not hyperbolic and contains quasi-flats if the genus of $S$ is 
greater than or equal to $4$, or the genus of $S$ is greater than or equal to $2$ and $\chi(S) \leq -6$.

\section{Genus ${\bf 2}$} \label{computation}

\begin{figure}[p]
\vspace{2 mm}
\centerline{\epsfxsize=3in \epsfbox{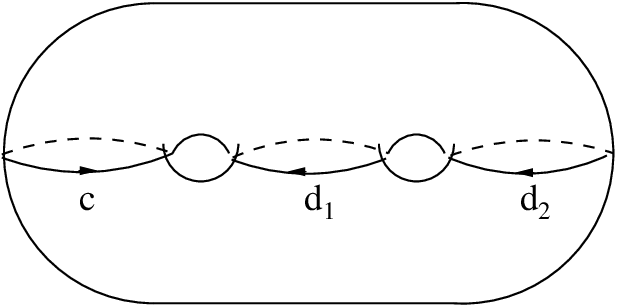}}
\caption{\sl Underlying curves $c$ and $d_1 \cup d_2$ in a genus $2$ surface (all weights are $1$)}
\label{genus2}
\vspace{2 mm}
\end{figure}

We consider the example of a closed orientable surface $S$ of genus $2$.  
A non trivial primitive homology class $\alpha$ 
can always be represented by a non separating simple closed curve with weight $1$.  Moreover,
a Seifert curve in a closed orientable surface of genus $2$, since its underlying curve is non separating,
can have at most two components.  Figure \ref{genus2} depicts multi-curves $c$ and $d_1 \cup d_2$ such
that $[[c]] = [[d_1]] + [[d_2]]$.

\begin{figure}[p]
\vspace{2 mm}
\centerline{\epsfxsize=3in \epsfbox{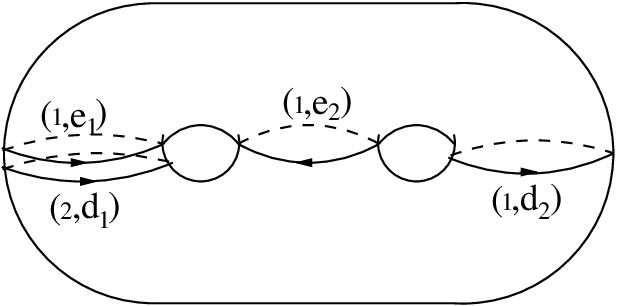}}
\caption{\sl Seifert curves $(2, 1, d)$ and $(1, 1, e)$}
\label{genus2kak}
\vspace{2 mm}
\end{figure}

We refer to a Seifert curve with one component as {\em type $1$} and a Seifert curve with
two components as {\em type $2$}.   Since $\alpha$ is primitive, a Seifert curve of type $1$ 
must have weight $1$.  It follows that distinct Seifert curves of type $1$ must intersect.
A Seifert curve of type $1$ can be disjoint from  a Seifert curve of type $2$, see Figure \ref{genus2} and 
distinct Seifert curves of type $2$ can be disjoint, see Figure \ref{genus2kak}.

Let $c$ be the underlying curve of a Seifert curve of type $1$ and $d = d_1 \cup d_2$ the underlying
curve of a Seifert curve
of type $2$ that are disjoint.   Then the three disjoint simple closed curves $c \cup d$ cut $S$ into pairs of 
pants.  Any Seifert curve that is disjoint from $c \cup d$ must have underlying curve
parallel to either $c$ or $d$.  Note that, since the weight of $c$ is $1$,
the weights for $d_1, d_2$ must also be $1$.  

Consider
the link of $[(1,c)]$ in $Kak(S, \alpha)$.  It consists of equivalence classes of Seifert curves of type $2$.
The Seifert curves of type $2$ have underlying curves that are
pairs of curves 
lying in $S \backslash c$, aren't parallel to $c$, and
are separating in $S \backslash c$ but not in $S$.  
There are infinitely many such pairs of curves.
More specifically, $S \backslash c$ is a twice punctured torus, 
so the curves are parallel curves that separate the two punctures and can be parametrized by ${\bf Q}$.
Distinct such curves can't be isotoped to be disjoint and hence correspond to 
distance two vertices of $Kak(S, \alpha)$.
This confirms that $Kak(S, \alpha)$ has dimension $1 = (2)(2) - 3$ near $[(1, c)]$, 
as prescribed by Theorem \ref{closedkakdim}.  


\begin{figure}[p]
\vspace{2 mm}
\centerline{\epsfxsize=3in \epsfbox{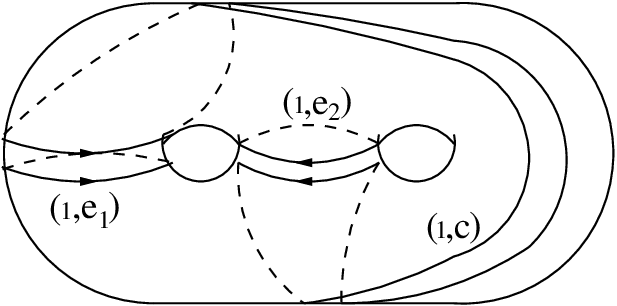}}
\caption{\sl Seifert curves $(1, 1, e)$ and $(1, c)$}
\label{genus2kakedge}
\vspace{2 mm}
\end{figure}

For $d = d_1 \cup d_2$ as in Figure \ref{genus2} or \ref{genus2kak}, 
we consider $S \backslash (d_1 \cup d_2)$,
a sphere with four punctures.  
The link of $[(w_1, w_2, d_1\cup d_2)]$ contains isotopy classes of Seifert curves of type $1$.
These are essential curves
that are separating in $S \backslash (d_1 \cup d_2)$ but not in $S$ and 
that partition the punctures of $S \backslash (d_1 \cup d_2)$ appropriately.  
There are infinitely many such curves.  They too can be parametrized by ${\bf Q}$.  
Note that distinct Seifert curves of type $1$ can't be isotoped to be disjoint and hence correspond to vertices of $Kak(S, \alpha)$ of distance two or more.  

In addition, the link of $[(w_1, w_2, d_1\cup d_2)]$ contains vertices $[(u_1, u_2, e_1\cup e_2)]$
such that one component of $e_1 \cup e_2$, say $e_1$, is parallel to a component of $d_1 \cup d_2$,
say $d_1$ and  $S \backslash (d_1 \cup d_2 \cup e_1 \cup e_2)$ consists of two pairs of pants 
and one annulus.   Seifert curves of this type can also be parametrized by $Q$, since $e_2$ is a curve
in a twice punctured torus that partitions the punctures appropriately and $e_1$ is parallel to $d_1$.
Note that, since the weights of $d_1, d_2$ are $w_1, w_2$, we must have 

\[\begin{array}{cccccccc}
w_1 & = & u_1 &  \pm & u_2 \\
w_2 &= & u_2 & & 
\end{array}\]

In summary, $Kak(S, \alpha)$ is a tree each of whose vertices has a countably infinite discrete 
({\it i.e.}, 0-dimensional) link.  

Recall that Johnson, Pelayo and Wilson showed that the Kakimizu complex of a knot in the 3-sphere 
is quasi-Euclidean.  The Kakimizu complex of the genus 2 surface
is an infinite graph, thus  Gromov hyperbolic.  In particular, it is not quasi-Euclidean.

\section{3-manifolds} \label{3d}

The definitions given for Seifert curve, infinite cyclic cover, Kakimizu complex and so forth carry over to codimension $1$ 
submanifolds in manifolds of any dimension.  In particular, they carry over to Seifert surfaces and Kakimizu complexes in the 
context of compact (possibly closed) $3$-manifolds.  One need merely replace ones by twos and twos by threes.  Instead of Seifert
curves, one considers Seifert surfaces.  Seifert surfaces are weighted essential surfaces that represent a given 
relative second homology class and have connected complement.  This ties into and generalizes some of the work in \cite{PS}.  

Let $S$ be a compact oriented surface.  Take $M = S \times I$.  Incompressible surfaces in a product manifold are either horizontal or vertical.   
Vertical surfaces have the form $c \times I$, where $c$ is a multi-curve in $S$.  
It follows that $Kak(M, [[c \times I]]) = Kak(S, [[c]])$, where $[[\cdot ]]$ denotes the homology class of $\cdot$.  

\begin{thm}
There exist $3$-manifolds with Gromov hyperbolic Kakimizu complex.  
\end{thm}

\proof
Let $S$ be the closed oriented surface of genus $2$, $\alpha$ a primitive homology class in $H_1(S)$ and
$c$ a compact $1$-manifold representating $\alpha$.  
Then $Kak(S, \alpha)$ is the graph discussed in Section \ref{computation}.  
In particular, $Kak(S, \alpha)$ is quasi-hyperbolic.  Take $M = S \times I$.   Then $Kak(M, [[c \times I]]) = Kak(S, \alpha)$ is
also quasi-hyperbolic.
\qed


\vspace{2 mm}

\noindent
Department of Mathematics

\noindent
1 Shields Avenue

\noindent
University of
California, Davis

\noindent
Davis, CA 95616

\noindent
USA

\end{document}